\documentclass[a4paper,10pt]{article}
\usepackage{amsmath,amsthm,amssymb,amsfonts}
\usepackage{enumerate,amscd,amsxtra,MnSymbol}
\usepackage{mathrsfs}
\usepackage{color}
\usepackage[cmtip,all]{xy}

\newtheorem{theorem}{Theorem}[section]
\newtheorem*{theorem*}{Theorem}
\newtheorem{lemma}[theorem]{Lemma}
\newtheorem{proposition}[theorem]{Proposition}
\newtheorem*{proposition*}{Proposition}
\newtheorem{corollary}[theorem]{Corollary}
\newtheorem*{corollary*}{Corollary}

\newtheorem{definition}[theorem]{Definition}
\newtheorem{definition*}{Definition}
\newtheorem{remark}[theorem]{Remark}
\newtheorem{remark*}{Remark}

\newtheorem{conjecture*}{Conjecture}

\newcommand{\caE}{{\mathcal E}}

\newcommand{\caH}{{\mathcal H}}

\newcommand{\caP}{{\mathcal P}}

\newcommand{\caG}{{\mathcal G}}

\newcommand{\naturals}{{\mathbb N}}
\newcommand{\integers}{{\mathbb Z}}

\newcommand{\uHom}{\underline{\mathrm{Hom}}}

\newcommand{\op}{\mathrm{op}}

\newcommand{\Hom}{\mathrm{Hom}}

\newcommand{\Sp}{\mathrm{Sp}}

\newcommand{\hocolim}{\mathrm{hocolim}}

\newcommand{\GW}{\mathrm{GW}}

\newcommand{\Cpx}{\mathrm{Cpx}}

\newcommand{\id}{\mathrm{id}}

\newcommand{\Cat}{\mathsf{Cat}}

\newcommand{\Set}{\mathsf{Set}}
\newcommand{\sSet}{\mathsf{sSet}}

\newcommand{\map}{\mathrm{map}}

\newcommand{\Pic}{\mathrm{Pic}}

\newcommand{\Spc}{\mathsf{Spc}}

\newcommand{\Mon}{\mathrm{Mon}}

\newcommand{\Fun}{\mathrm{Fun}}
\newcommand{\SymMonCat}{\mathsf{SymMonCat}}

\newcommand{\preadd}{\mathrm{preadd}}

\newcommand{\Fin}{\mathsf{Fin}}

\newcommand{\lax}{\mathrm{lax}}

\newcommand{\Perf}{\mathrm{Perf}}
\newcommand{\st}{\mathrm{st}}
\newcommand{\Ar}{\mathrm{Ar}}
\renewcommand{\CD}{\mathrm{CD}}
\newcommand{\str}{\mathrm{str}}
\newcommand{\Tw}{\mathrm{Tw}}
\newcommand{\Hyp}{\mathrm{Hyp}}

\newcommand{\lex}{\mathrm{lex}}
\newcommand{\Span}{\mathrm{Span}}
\newcommand{\KR}{\mathrm{KR}}

\title{A Grothendieck-Witt space for stable infinity categories with duality}
\author{Markus Spitzweck}

\begin{document}

\maketitle

\begin{abstract}

We construct a Grothendieck-Witt space for any stable infinity category with duality.
If we apply our construction to perfect complexes over a commutative ring in which $2$ is invertible
we recover the classical Grothendieck-Witt space. Our Grothendieck-Witt space is a grouplike $E_\infty$-space
which is part of a genuine $C_2$-spectrum, the connective real $K$-theory spectrum.

\end{abstract}

\tableofcontents

\section{Introduction}

In this paper we carry over the hermitian $S_\bullet$-construction which can be found e.g. in
\cite{schlichting.mv} to the $\infty$-categorical setting. The input of our construction
is an $\infty$-category with duality in the sense of \cite{hls} whose underlying $\infty$-category
is stable. For such an $\infty$-category $C$ we define a Grothendieck-Witt space $\GW(C)$
which has the structure of a grouplike $E_\infty$-space, so gives rise to a connective spectrum.
We show in the last section that this spectrum is in fact part of a genuine $C_2$-spectrum
$\KR(C)$, the connective real $K$-theory spectrum of $C$.

This way we obtain for example for any $E_\infty$-ring spectrum $R$ and tensor invertible object
$L \in \Perf(R)$ (here $\Perf(R)$ denotes the stable $\infty$-category of perfect $R$-modules)
a real $K$-theory spectrum $\KR(R,L)$, in particular spectra $\KR(R,R[n])$ for any $n \in \integers$,
by considering the $L$-twisted duality on $\Perf(R)$ (see \cite[\S 8]{hls}).

To justify our constructions we prove in section \ref{j53454t} that our Grothendieck-Witt space
is equivalent to the classical Grothendieck-Witt space (as defined e.g. in \cite{schlichting.mv})
in the case that $R$ is a discrete ring in which $2$ is invertible and $L$ a shifted invertible
(in the discrete sense) $R$-module.

\vskip.7cm

{\bf Acknowledgements:}
I would like to thank Hongyi Chu,
David Gepner, Ha\-drian Heine, Kristian Moi, Thomas Nikolaus, Oliver R\"ondigs, Manfred Stelzer, Sean Tilson and Girja Tripathi for very helpful
discussions and suggestions on the subject.

\section{Recollections and preliminaries}

We use the same conventions as in \cite{hls}. In particular $\Cat_\infty^{hC_2}$ is the $\infty$-category
of small $\infty$-categories with duality.
We denote be $\Cat_\infty^\st$ the subcategory of $\Cat_\infty$ of stable $\infty$-categories and
exact functors between them.

It follows from \cite[Proposition 2.2]{hls} that the induced functor
$$(\Cat_\infty^\st)^{hC_2} \to \Cat_\infty^{hC_2}$$
(on the left we use the induced $C_2$-action) is a monomorphism in $\widehat{\Cat}_\infty$
whose good image consists of those $\infty$-categories with duality whose underlying
$\infty$-category is stable and those functors between $\infty$-categories with duality whose
underyling functor is exact. We write for this good image $(\Cat_\infty^{hC_2})^\st$.

Usually we will not distinguish between a category and its nerve viewed as an $\infty$-category.

We will frequently see objects $[n] \in\Delta$ as categories.
For a category $C$ we write $\Ar(C)$ for the arrow category, i.e. the functor category
$\Fun([1],C)$.

\begin{proposition}
\label{hte4z44}
Let $n \in \naturals$ and $C_1,\ldots,C_n$ be stable $\infty$-categories. Let $D$ be an $\infty$-category
which admits finite limits and denote by $\Sp(D)$ the stabilization of $D$. Let
$$\Fun'(C_1 \times \cdots \times C_n,D) \subset \Fun(C_1 \times \cdots \times C_n,D)$$
be the full subcategory on those functors which preserve finite limits separately in each variable, and let
$$\Fun'(C_1 \times \cdots \times C_n,\Sp(D)) \subset \Fun(C_1 \times \cdots \times C_n,\Sp(D))$$
be the full subcategory on those functors which are exact separately in each variable.
Then composition with the functor $\Omega^\infty \colon \Sp(D) \to D$ induces an equivalence
$$\Fun'(C_1 \times \cdots \times C_n,\Sp(D)) \to \Fun'(C_1 \times \cdots \times C_n,D)$$
of $\infty$-categories. If $C_1= \cdots =C_n$ then this equivalence respects the $\Sigma_n$-actions.
\end{proposition}

\begin{proof}
This follows from \cite[Corollary 1.4.2.23.]{lurie.higheralgebra}.
\end{proof}

Let $C \in \Cat_\infty$. To give a duality on $C$ (or equivalently on $C^\op$)
is the same as to give a $C_2$-homotopy fixed point of $\Fun(C \times C,\Spc)$
(or equivalently of $\Fun(C^\op \times C^\op,\Spc)$) which is underlying a perfect pairing
(see \cite[Corollary 7.3]{hls}). If $\varphi \in \Fun(C^\op \times C^\op,\Spc)$ corresponds to
a duality then $\varphi$ is informally given by $(X,Y) \mapsto \map(X,Y^\vee)$.

If $C$ is now stable it follows from Proposition \ref{hte4z44} that such a $C_2$-homotopy fixed point
is the same as a $C_2$-homotopy fixed point $B$ of $\Fun(C^\op \times C^\op, \Sp)$ which is
nondegenerate representable in the sense of \cite{lurie.L4}.
The bilinear functor $B$ is informally given by $(X,Y) \mapsto \map^\Sp(X,Y^\vee)$,
where $\map^\Sp$ denotes the mapping spectrum functor.

\section{The Grothendieck-Witt space}
\label{jtez4t4x}

Let $C \in \Cat_\infty^\st$. Building on Waldhausen's definition we define
for any $n \in \naturals$ the $\infty$-category $S_n(C)$ to be the full subcategory of
the functor category $\Fun(\Ar([n]),C)$ on those functors $A$ such that for any $0 \le i \le n$
the object $A_{i,i}$ is a zero object in $C$ and such that for any $0 \le i \le j \le k \le n$
the square
$$\xymatrix{A_{i,j} \ar[r] \ar[d] & A_{i,k} \ar[d] \\
A_{j,j} \ar[r] & A_{j,k}}$$
is exact in $C$.
These properties are preserved by the suspension and loop functors on $C$, thus the $\infty$-categories
$S_n(C)$ are stable. The simplicial $\infty$-category $\Fun(\Ar([\bullet]),C)$
restricts to a simplicial $\infty$-category $S_\bullet(C)$. Taking levelwise core groupoids yields
the simplicial object $S_\bullet^\sim(C)$ in spaces whose realization we denote by $|S_\bullet^\sim(C)|$.

\begin{definition}
The $K$-theory space $K(C)$ of the stable $\infty$-category
$C$ is defined to be the loop space $\Omega |S_\bullet^\sim(C)|$, where we take
a zero object of $C$ as base point.
\end{definition}

\begin{remark}
The space $K(C)$ has the natural structure of a grouplike $E_\infty$-space (see also the discussion at the end of this section
for the case of the Grothendieck-Witt space).
\end{remark}

Our $\infty$-categorical definition of the Grothendieck-Witt space of a stable $\infty$-category
with duality is modelled on the hermitian $S_\bullet$-construction given for example in
\cite{schlichting.mv}. This uses the edgewise subdivision of a simplicial object which we introduce now.

\begin{definition}
Let $X \colon \Delta^\op \to C$, $[n] \mapsto X_n$, be a simplicial object in an $\infty$-category $C$.
Then the edgewise subdivision $E(X)$ is defined to be the simplicial object $X \circ \iota^\op$, where
$\iota \colon \Delta \to \Delta$ is the endofunctor defined by $[n] \mapsto [n]^\op * [n]$.
\end{definition}

Thus we have $E(X)_n = X_{2n+1}$.
The inclusions $[n] \hookrightarrow [n]^\op * [n]$ define a natural transformation from the identity functor
on $\Delta$ to $\iota$ and thus we are at the disposal of a natural map of simplicial objects $E(X) \to X$.

For $C$ a stable $\infty$-category we let $S_\bullet^e(C):=E(S_\bullet(C))$, and likewise
$S_\bullet^{e,\sim}(C):=E(S_\bullet^\sim(C))$.

For each $n \in \naturals$ the category $[n]$ has a unique structure of a category with stict duality and
the assignment $[n] \mapsto [n]^\op * [n]$ can be viewed as a functor from
$\Delta$ to the category of
categories with strict duality $\CD$, therefore the same holds for the assignment
$[n] \mapsto \Ar([n]^\op * [n])$.

In \cite[\S 11]{hls} a functor $\epsilon \colon \CD \to \Cat_\infty^{hC_2}$ is constructed.
Moreover $\Cat_\infty^{hC_2}$ is cartesian closed, and the internal hom commutes with the forgetful
functor $\Cat_\infty^{hC_2} \to \Cat_\infty$.

Thus for an $\infty$-category $C$ with duality $\Fun(\Ar([n]^\op * [n]),C)$ is an object of $\Cat_\infty^{hC_2}$ functorial in $[n]$.

If now $C$ is a stable $\infty$-category
with duality then for any $n \in \naturals$ the full subcategory $S_n(C)$ of $\Fun(\Ar([n]),C)$ is preserved
by the duality (since the dual of an exact square is again
an exact square), thus $S_\bullet^e(C)$ can be viewed as a simplicial object in $\Cat_\infty^{hC_2}$,
and $S_\bullet^{e,\sim}(C)$ can be viewed as a simplicial object in $\Spc^{hC_2} \simeq \Spc[C_2]$.

Taking levelwise the homotopy $C_2$-fixed points of the latter object defines the simplicial
space $(S_\bullet^{e,\sim}(C))_h$.

\begin{definition}
For $C$ a stable $\infty$-category with duality we let the Grothendieck-Witt space $\GW(C)$ of $C$ be the homotopy fiber of the
composition $$|(S_\bullet^{e,\sim}(C))_h| \to |S_\bullet^{e,\sim}(C)| \to |S_\bullet^\sim(C)|.$$
\end{definition}

We now equip $\GW(C)$ with an $E_\infty$-structure which will turn out to be grouplike (i.e. an infinite loop space
structure), see Proposition \ref{wet5u6}.
The monomorphisms $$\Cat_\infty^\st \to \Cat_\infty^\preadd \to \Cat_\infty$$
as well as the full embedding $$\Cat_\infty^\preadd \to \SymMonCat_\infty$$
of $\infty$-categories carry $C_2$-actions (see \cite[\S 6]{hls}), hence we have an induced composition
$$s \colon (\Cat_\infty^{hC_2})^\st \simeq (\Cat_\infty^\st)^{hC_2} \to (\Cat_\infty^\preadd)^{hC_2} \to \SymMonCat_\infty^{hC_2}.$$

Thus, since $S_\bullet^e(C)$ is in fact naturally a simplicial object in $(\Cat_\infty^{hC_2})^\st$,
we obtain a simplicial object $s(S_\bullet^e(C))$ in $\SymMonCat_\infty^{hC_2}$,
and applying the functor $$\SymMonCat_\infty^{hC_2} \overset{(-)^\sim}{\longrightarrow} \Mon_{E_\infty}(\Spc)[C_2]$$
yields a lift of $S_\bullet^{e,\sim}(C)$ to a simplicial object of $\Mon_{E_\infty}(\Spc)[C_2]$.

Similarly the map $$S_\bullet^{e,\sim}(C) \to S_\bullet^\sim(C)$$ lifts to a map between simplicial objects
in $\Mon_{E_\infty}(\Spc)$.

Denoting the lifts with the same symbols we obtain maps
$$(S_\bullet^{e,\sim}(C))_h \to S_\bullet^{e,\sim}(C) \to S_\bullet^\sim(C)$$
of simplicial objects in $\Mon_{E_\infty}(\Spc)$. Taking realizations and the fiber of the induced composition equips
$\GW(C)$ with a natural $E_\infty$-structure.

\section{The comparison}
\label{j53454t}

We denote the (hermitian) $S_\bullet$-construction used in \cite{schlichting.mv} for an exact category with weak equivalences
(and duality) $\caE$ by the same symbols as we used in the $\infty$-categorical situation except that we write
$S^\str$ instead of $S$. Thus for example if $\caE$ has a duality then $S_\bullet^{\str,e}(\caE)$ is a simplicial
exact category with weak equivalences and duality and $S_\bullet^{\str,e,\sim}(\caE)$ denotes the simplicial
subcategory of weak equivalences.

We denote by $$|\_| \colon \Cat^1 \to \Spc$$ the natural functor form the $1$-category of small categories
$\Cat^1$ to the $\infty$-category of spaces which takes the realization of the nerve.
Thus if $C_\bullet$ is for example a simplicial category then $|C_\bullet|$ will be a simplicial object in $\Spc$. 
The realization of this simplicial object is denoted by $|C_\bullet|_r$.

The Grothendieck-Witt space $\GW(\caE)$ is then defined to be the homotopy fiber of the natural map
$$|(S_\bullet^{\str,e,\sim}(\caE))_h|_r \to |S_\bullet^{\str,\sim}(\caE)|_r.$$
As in the $\infty$-categorical case we can equip $\GW(\caE)$ with a natural $E_\infty$-structure
(use that the functor $\caH^\lax$ (see \cite[\S 11]{hls}) is symmetric monoidal for the cartesian symmetric monoidal
structures since it is a right adjoint).

For a ring $R$ we denote by $\caP_R$ the category of finitely generated projective $R$-modules and
by $\Cpx^b(\caP_R)$ the exact category with weak equivalences of bounded complexes with values in $\caP_R$.
We denote by $\Perf(R)$ the stable $\infty$-category of perfect $R$-modules.

Note that we exhibit a natural functor $\Cpx^b(\caP_R) \to \Perf(R)$ which is a localization
at the quasi isomorphisms.

We now assume that $R$ is commutative, fix for the whole section an integer $N \in \integers$
and an invertible $R$-module $L$ and equip $\Cpx^b(\caP_R)$ with the strong duality
$X \mapsto \uHom(X,L[N])$.

\cite[Corollary 8.5]{hls} equips $\Perf(R)$ with the duality given by the object $L[N] \in \Pic(\Perf(R))$,
and the naturality of the construction of loc. cit. shows that the functor
$$\Cpx^b(\caP_R) \to \Perf(R)$$
preserves the dualities.

To emphsize the dependence on the duality we denote the corresponding
Grothendieck-Witt spaces by $\GW(\Cpx^b(\caP_R),N,L)$ and $\GW(\Perf(R),N,L)$.

For any $n \in \naturals$ we obtain a functor $$\Fun(\Ar([n]),\Cpx^b(\caP_R)) \to
\Fun(\Ar([n]),\Perf(R))$$
between $\infty$-categories with duality.

The restriction to the full subcategory $S_n^\str(\Cpx^b(\caP_R))$ of this functor factors through
$S_n(\Perf(R))$ yielding functors
$$S_n^\str(\Cpx^b(\caP_R)) \to S_n(\Perf(R))$$
\begin{equation}
\label{sdrgt}
S_n^{\str,\sim}(\Cpx^b(\caP_R)) \to S_n^\sim(\Perf(R))
\end{equation}
between $\infty$-categories with duality.

The functors $$\caH^\lax(S_n^{\str,\sim}(\Cpx^b(\caP_R))) \to \caH^\lax(S_n^\sim(\Perf(R)))$$
induced on lax hermitian objects (see \cite[\S 11]{hls}) by the latter functors
together with the equivalence $$\caH(S_n^\sim(\Perf(R))) \simeq \caH^\lax(S_n^\sim(\Perf(R)))$$ yields functors
\begin{equation}
\label{h54wergrh}
(S_n^{\str,\sim}(\Cpx^b(\caP_R)))_h \to (S_n^\sim(\Perf(R)))_h
\end{equation}
(see also \cite[Proposition 11.8]{hls}).

Every map in $S_n^{\str,\sim}(\Cpx^b(\caP_R))$ is sent to an equivalence under the functor (\ref{sdrgt}),
thus we obtain maps
\begin{equation} 
\label{htr46z}
|S_n^{\str,\sim}(\Cpx^b(\caP_R))| \to S_n^\sim(\Perf(R))
\end{equation}
in $\Spc[C_2]$. Also every map in $(S_n^{\str,\sim}(\Cpx^b(\caP_R)))_h$ is sent to an equivalence
under the functor (\ref{h54wergrh}), thus we obtain
maps
\begin{equation}
\label{htr3eetg}
|(S_n^{\str,\sim}(\Cpx^b(\caP_R)))_h| \to (S_n^\sim(\Perf(R)))_h
\end{equation}
in $\Spc$.

After edgewise subdivision we get a map
$$|(S_\bullet^{\str,e,\sim}(\Cpx^b(\caP_R)))_h| \to (S_\bullet^{e,\sim}(\Perf(R)))_h$$
between simplicial objects in $\Spc$.

Altogether we arrive at a 
commutative square
$$\xymatrix{|(S_\bullet^{\str,e,\sim}(\Cpx^b(\caP_R)))_h| \ar[r] \ar[d] & 
|S_\bullet^{\str,\sim}(\Cpx^b(\caP_R))| \ar[d] \\
(S_\bullet^{e,\sim}(\Perf(R)))_h \ar[r] & S_\bullet^\sim(\Perf(R))}$$
of simplicial objects in $\Spc$, which yields after taking realizations
and fibers of the horizontal induced maps the
comparison map
\begin{equation}
\label{ht44rtg}
\GW(\Cpx^b(\caP_R),N,L) \to \GW(\Perf(R),N,L).
\end{equation}

\begin{remark}
The comparison map (\ref{ht44rtg}) can be made compatible with the $E_\infty$-structures on both sides.
We leave the details to the interested reader.
\end{remark}

\begin{theorem}
If $2$ is invertible in $R$ the comparison map (\ref{ht44rtg}) is an equivalence.
\end{theorem}

\begin{proof}
Combine the next two Lemmas.
\end{proof}

\begin{lemma}
\label{hte44rt}
The maps (\ref{htr46z}) are equivalences.
\end{lemma}
\begin{proof}
This is standard.
\end{proof}

The main input to our comparison statement is

\begin{lemma}
\label{hrergr}
If $2$ is invertible in $R$ then the maps (\ref{htr3eetg})
are equivalences.
\end{lemma}

\begin{proof}
We have a commutative diagram
\begin{equation}
\label{grert344}
\xymatrix{|(S_n^{\str,\sim}(\Cpx^b(\caP_R)))_h| \ar[r] \ar[d] &
S_n^\sim(\Perf(R))_h \ar[d] \\
|S_n^{\str,\sim}(\Cpx^b(\caP_R))| \ar[r] &
S_n^\sim(\Perf(R))}
\end{equation}
in $\Spc$.
We want to show that the upper horizontal map is an equivalence.
By Lemma \ref{hte44rt} the lower horizontal map is an equivalence.
We will show that for any $X \in S_n^{\str,\sim}(\Cpx^b(\caP_R))$
the space of paths $\map(X,X^\vee)$ in $|S_n^{\str,\sim}(\Cpx^b(\caP_R))|$ (or equivalently in
$S_n^\sim(\Perf(R))$) carries a natural $C_2$-action, that the homotopy fibers of the vertical
maps in the diagram over $X$ (resp. the image of $X$) are canonically identified with
$\map(X,X^\vee)^{hC_2}$ and that the induced map (by the commutative square) on these fibers respect
these identifications. From this the claim follows.

We first apply \cite[Lemma 4]{schlichting.mv} to the exact category $S_n^\str(\Cpx^b(\caP_R))$
with duality (considering only the isomorphisms
as weak equivalences) to obtain a category $C=S_n^\str(\Cpx^b(\caP_R))_{\mathrm{iso}}^\str$ with a strict duality
which is equivalent to $S_n^\str(\Cpx^b(\caP_R))$ as category with (strong) duality.
We denote by $C^\sim$ the subcategory of $C$ of weak equivalences (which correspond to the objectwise quasi isomorphisms
in $S_n^\str(\Cpx^b(\caP_R))$).

For a category $D$ we let $\Tw(D)$ be the twisted arrow category of $D$ whose objects are the morphisms of $D$,
and a map from $f \colon A \to B$ to $g \colon C \to D$ is a commutative square
$$\xymatrix{A \ar[r] \ar[d]^f & C \ar[d]^g \\
B & D \ar[l]}$$
in $D$.
If $D$ has a strict duality then the assignment $f \mapsto f^\vee$ defines a (strict) $C_2$-action
on $\Tw(D)$ whose (strict) $C_2$-fixed points is the category of hermitian objects of $D$.

Similarly for an $\infty$-category $D$ the assignment $$\Delta^\op \ni [n] \mapsto \map([n] * [n]^\op,D)$$
defines a complete Segal space whose associated $\infty$-category is defined to be the twisted
arrow category $\Tw(D)$ of $D$
(this is compatible with the $1$-categorical definition). If $D$ has a duality then the above assignment
has values in $\Spc[C_2]$, thus $\Tw(D)$
has a $C_2$-action. Moreover by the construction in \cite[\S 11]{hls}
we have a canonical equivalence $$\Tw(D)^{hC_2} \simeq \caH^\lax(D).$$
The canonical map $\Tw(D) \to D \times D^\op$ is $C_2$-equivariant, where the action
on $D \times D^\op$ is given by $(X,Y) \mapsto (Y^\vee,X^\vee)$,
and we have $(D \times D^\op)^{hC_2} \simeq D$.

The right vertical map of diagram (\ref{grert344}) can thus be identified with the map
$$\Tw(S_n^\sim(\Perf(R)))^{hC_2} \to (S_n^\sim(\Perf(R)) \times S_n^\sim(\Perf(R))^\op)^{hC_2}.$$
We have a commutative diagram
$$\xymatrix{\Tw(C^\sim) \ar[r] \ar[d] & \Tw(S_n^\sim(\Perf(R))) \ar[d] \\
C^\sim \times (C^\sim)^\op \ar[r] & S_n^\sim(\Perf(R)) \times S_n^\sim(\Perf(R))^\op}$$
in $\Cat_\infty[C_2]$.
In the induced diagram
$$\xymatrix{|\Tw(C^\sim)| \ar[r] \ar[d] & \Tw(S_n^\sim(\Perf(R))) \ar[d] \\
|C^\sim \times (C^\sim)^\op| \ar[r] & S_n^\sim(\Perf(R)) \times S_n^\sim(\Perf(R))^\op}$$
in $\Spc[C_2]$
the horizontal maps are equivalences.
Thus diagram (\ref{grert344}) can be identified with the diagram
\begin{equation}
\label{ht4436z}
\xymatrix{|\Tw(C^\sim)^{C_2}| \ar[r] \ar[d] & |\Tw(C^\sim)|^{hC_2} \ar[d] \\
|(C^\sim \times (C^\sim)^\op)^{C_2}| \ar[r] & |C^\sim \times (C^\sim)^\op|^{hC_2}}
\end{equation}
whose lower entries can be identified with $|C^\sim|$.

For $X \in C$ we let $P_X$ be defined by the (strict) pullback diagram
\begin{equation}
\label{j64236sf}
\xymatrix{P_X \ar[r] \ar[d] & \Tw(C^\sim) \ar[d] \\
(C^\sim \times (C^\sim)^\op)/(X,X^\vee) \ar[r] & C^\sim \times (C^\sim)^\op}
\end{equation}
of categories.
Since $(X,X^\vee) \in C^\sim \times (C^\sim)^\op$ is a fixed point with respect to the $C_2$-action
$P_X$ inherits a $C_2$-action and this diagram becomes $C_2$-equivariant.
Taking $C_2$-fixed points of this diagram gives a diagram canonically isomorphic to the pullback diagram
\begin{equation}
\label{h2356u6eg}
\xymatrix{P_X^{C_2} \ar[r] \ar[d] & C^\sim_h \ar[d] \\
C^\sim/X \ar[r] & C^\sim.}
\end{equation}
For a map $X \to Y$ in $C^\sim$ we have an induced $C_2$-equivariant map $P_X \to P_Y$.

\vskip.2cm

{\it Claim 1:} For any map $X \to Y$ in $C^\sim$ the map $|P_X| \to |P_Y|$ is an equivalence.

\vskip.1cm

{\it Claim 2:} For any $X \in C$ the map $|P_X^{C_2}| \to |P_X|^{hC_2}$ is an equivalence.

\vskip.1cm

{\it Claim 3:} For any map $X \to Y$ in $C^\sim$ the map $|P_X^{C_2}| \to |P_Y^{C_2}|$ is an equivalence.

\vskip.2cm

Claim 3 follows from Claims 1 and 2.

It follows from Claim 1 and Quillen's Theorem B (dual of \cite[Theorem 5.6]{goerss-jardine}) that the realization
of diagram (\ref{j64236sf}) is a pullback diagram, similarly it follows from Claim 3 and Quillen's Theorem B
that the realization
of diagram (\ref{h2356u6eg}) is a pullback diagram.

Thus the induced map on the homotopy fibers over an $X \in C$ of the vertical maps in diagram (\ref{ht4436z}) can be identified with
with the map $|P_X^{C_2}| \to |P_X|^{hC_2}$ (for the second fiber note that
homotopy fixed points preserve fiber sequences)
which is an equivalence by Claim 2. So we see that if we prove Claims 1 and 2
the proof is finished.

{\it Proof of Claim 1:} It follows from \cite[Propositions 6.2 and 8.2]{dwyer-kan.calculating}
and the correction \cite{dugger.class} that $|P_X|$ is canonically
equivalent to the mapping space $\map(X,X^\vee)$ in $|C^\sim|$:
The category $P_X$ is naturally isomorphic to the category denoted $C^\sim(X,X^\vee)_{\Hom - \mathit{tw}}$ in \cite{dugger.class}),
and this a collection of connected components of $C(X,X^\vee)_{\Hom - \mathit{tw}}$.
The factorizations necessary for these arguments are given by cylinder constructions.
Thus the claim follows.

{\it Proof of Claim 2:} Let $X \in C$ and $X'$ be the image of $X$ in $S_n^\str(\Cpx^b(\caP_R))$.
Let $N \integers[\Delta^\bullet]$ be the cosimplicial
object in $\Cpx^b(\caP_\integers)$ which assigns to $[n]$ the complex corresponding
to the simplicial abelian group $\integers[\Delta^\bullet]$ under the Dold-Kan correspondence.
Thus it is Reedy cofibrant, and the cosimplicial object $N \integers[\Delta^\bullet] \otimes X'$
in $S_n^\str(\Cpx^b(\caP_R))$ is a special cosimplicial resolution of $X'$ in the sense of
\cite{dwyer-kan.function}.
We let $X^\bullet$ be the image of this cosimplicial object in $C$.
Also let $(C^\sim/X)_f$ be the full subcategory of $C^\sim/X$ on those maps $Y \to X$
which are surjections.
Then by \cite[Proposition 6.12]{dwyer-kan.function}
the functor $\varphi \colon \Delta \to (C^\sim/X)_f$ which sends $[n]$ to $X^n \twoheadrightarrow X$ is left cofinal.

Let the categories $R$ and $S$ be defined by the pullback diagram
\begin{equation}
\label{hre32z4t}
\xymatrix{S \ar[r]^\psi \ar[d] & R \ar[r]^i \ar[d] & P_X^{C_2} \ar[d] \\
\Delta \ar[r]^\varphi & (C^\sim/X)_f \ar[r] & C^\sim/X.}
\end{equation}
Let $r \in R$, so $r$ consists of a surjection $Y \twoheadrightarrow X$ in $C^\sim$ together with
a hermitian structure on $Y$. Since for any map $Z \to Y$ in $C^\sim$ there exists a unique hermitian
structure on $Z$ compatible with the one on $Y$ it follows that the natural functor
$$\psi/r \to \varphi/(Y \twoheadrightarrow X)$$ is an isomorphism, hence $\psi$ is also left cofinal.

The vertical functors in diagram (\ref{hre32z4t}) are right fibrations, and the fiber over an object
$Y \to X$ is $\Hom_{C^\sim}(Y,Y^\vee)^{C_2}$ (i.e. the set of hermitian structures on $Y$). Thus for a vertical map
$A \to B$ in this diagram we exhibit a functor $j_B \colon B^\op \to \Set \hookrightarrow \sSet$, and $|A|$ is naturally
equivalent to $\hocolim j_B$. It follows that $|S| \to |R|$ is an equivalence.

Let $K$ be the simplicial set defined by $[n] \mapsto \Hom_C(X^n,(X^n)^\vee)$
($K$ has then in fact the structure of a simplicial $R$-module) and $K^\sim$ the subsimplicial set
on those simplices which are in $C^\sim$. The simplicial set $K$ has a natural $C_2$-action
and $K^\sim$ is stable under this action.

It follows from the above considerations that $|S|$ is naturally equivalent to $(K^\sim)^{C_2}$

Note that the natural map $$K^{C_2} \to K^{hC_2}$$ is an weak homotopy equivalence since $2$ is invertible
in $R$, thus, since $K^\sim$ consists of certain connected components of $K$, the same follows for the map $$(K^\sim)^{C_2} \to (K^\sim)^{hC_2}.$$

For a map $f \colon Y \to Z$ in $S_n^\str(\Cpx^b(\caP_R))$ denote by $c(f)$ the construction 
\cite[1.5.5]{weibel.hom} applied to the map $f^\op$ in $S_n^\str(\Cpx^b(\caP_R))^\op$. We therefore
obtain a factorization $Y \to c(f) \to Z$ of $f$ into an inclusion which is a quasi isomorphism followed by a surjection,
and moreover there is a retraction $c(f) \to Y$ of the first map. This construction is functorial in $f$.
If $Y$ has a hermitian structure then the retraction induces a hermitian structure on $c(f)$.
These constructions can be tranported to $C$.

For an object $a \in P_X^{C_2}$ with underlying object $f \colon Y \to X$ in $C^\sim/X$ the above factorization applied
to $f$ yields an object $(c(f) \to X) \in (C^\sim/X)_f$ and also an object in $R$ (using the induced hermitian structure
on $c(f)$). This assignment defines a functor $p \colon P_X^{C_2} \to R$ together with natural transformations
$\id \to i \circ p$ and $\id \to p \circ i$. It follows that the realization of $i$ is an equivalence.

Hence we have seen that the natural map $|S| \to |P_X^{C_2}|$ in $\Spc$ is an equivalence, and that in
the commutative square
$$\xymatrix{(K^\sim)^{C_2} \ar[r] \ar[d] & |P_X^{C_2}| \ar[d] \\
(K^\sim)^{hC_2} \ar[r] & |P_X|^{hC_2}}$$
in $\Spc$
the upper horizontal and the left vertical maps are equivalences.
Also the lower horizontal map is an equivalence.
Hence Claim 2 and thus the Lemma are proved.
\end{proof}

\section{The zeroth Grothendieck-Witt group}

Let $C \in \Cat_\infty^{hC_2}$. Then the right fibration $$p \colon \Tw(C) \to C \times C^\op$$
inherits a $C_2$-action (see the proof of Lemma \ref{hrergr}).
Thus for $X \in C$ the space $\map(X,X^\vee)$ has a natural $C_2$-action,
since it arises as the homotopy fiber of $p$ over a homotopy fixed point for the $C_2$-action.

Using the equivalence $\Tw(C)^{hC_2} \simeq \caH^\lax(C)$ we see that the fiber
over $X$ of the functor $\caH^\lax(C) \to C$ is the $\infty$-groupoid $\map(X,X^\vee)^{hC_2}$,
so to give a lax hermitian structure on $X$ is the same as to give a $C_2$-homotopy
fixed point of $\map(X,X^\vee)$.

On the other hand the symmetric functor $C^\op \times C^\op \to \Spc$ corresponding to the duality on $C$
can be viewed as a map in $\widehat{\Cat}_\infty[C_2]$ (where the $C_2$-action on the source is the switch map
and on the target the trivial action), and taking homotopy fixed points yields a functor
$$C^\op \to \Spc[C_2]$$ which is informally given by $X \mapsto \map(X,X^\vee)$.
This way $\map(X,X^\vee)$ also inherits a $C_2$-action which can be seen to be naturally equivalent
to the action from above.

If $C$ is stable the same argument as above yields a functor $$C^\op \to \Sp[C_2]$$
which is informally given by $X \mapsto \map^\Sp(X,X^\vee)$. Composing with $$\Omega^\infty \colon \Sp \to \Spc$$
yields the functor above. Let $Q \colon C^\op \to \Sp$ be given by $X \mapsto \map^\Sp(X,X^\vee)^{hC_2}$.
We see that a lax hermitian structure on an $X \in C$ is the same as a quadratic object structure of $(C,Q)$ on $X$
in the sense of \cite{lurie.L5}, and a hermitian structure on $X$ is same as a Poincare object structure on $X$.

\begin{lemma}
\label{rergrer}
Let $I \in \Cat_\infty$ and $C \in \Cat_\infty^{hC_2}$.
Then there is a natural functor
$$\caH(\Fun(\caG^\lax(I),C)) \to \Fun(I,\caH^\lax(C))$$
inducing an equivalence on core groupoids.
\end{lemma}

\begin{proof}
Functorially in $J \in \Cat_\infty$ we have a chain of maps 
$$\map(J,\caH(\Fun(\caG^\lax(I),C))) \simeq \map(\caG(J),\Fun(\caG^\lax(I),C))$$
$$\simeq \map(\caG(J) \times \caG^\lax(I),C) \to \map(\caG^\lax(J) \times \caG^\lax(I),C)$$
$$\to \map(\caG^\lax(J \times I),C) \simeq \map(J \times I,\caH^\lax(C)) \simeq
\map(J,\Fun(I,\caH^\lax(C)))$$
in $\Spc$ defining the functor in question.
Taking core groupoids the functor reduces to the equivalence
$$\map(\caG^\lax(I),C) \simeq \map(I,\caH^\lax(C)).$$
\end{proof}

\begin{remark}
In general the functor in Lemma \ref{rergrer} is not an equivalence, since in general
the map $\caG^\lax(J \times I) \to \caG(J) \times \caG^\lax(I)$ is not an equivalence.
We always have an equivalence $$\caH(\Fun(\caG(I),C)) \simeq \Fun(I,\caH(C))$$
of $\infty$-categories.
\end{remark}

Let $C \in \Cat_\infty^{hC_2}$ and $\alpha \in \caH(\Fun([1],C))$ (so $\alpha$ can be identified
with a lax hermitian object of $C$). Let $F_\alpha$ be fiber over $\alpha$ of the functor
$$\caH(\Fun([3],C)) \to \caH(\Fun([1],C))$$
induced by the duality preserving functor $[1] \to [3]$
which sends $0$ to $1$ and $1$ to $2$.

Because of Lemma \ref{rergrer} the core groupoid $F_\alpha^\sim$ can then be identified
with $(\caH^\lax(C)_{/\alpha})^\sim$ (use $\caG^\lax([0]) \simeq [1]$ and
$\caG^\lax([1]) \simeq [3]$). Since $\caH^\lax(C) \to C$ is a right fibration
the latter category can be identified with $(C_{/\alpha(0)})^\sim$.

The duality preserving functor $[3] \to [2]$ given by $0 \mapsto 0, 1 \mapsto 1, 2 \mapsto 1, 3 \mapsto 2$
exhibits $\Fun([2],C)$ as the full subcategory of $\Fun([3],C)$ of those functors for which the middle
induced map is an equivalence. Therefore we also get a full embedding
$$\caH(\Fun([2],C)) \to \caH(\Fun([3],C)).$$
The considerations above show that the core groupoid of the fiber of
the functor $$\caH(\Fun([2],C)) \to \caH(C)$$
induced by $[0] \to [2]$, $0 \mapsto 1$, over a hermitian object $X$
is naturally equivalent to $(C_{/X})^\sim$ (here $X$ also denotes the object underlying the hermitian object).

\begin{proposition}
Let $C \in (\Cat_\infty^{hC_2})^\st$ and $X \in \caH(C)$.
Let $F$ be the fiber over $X$ of the functor $$\caH(S_2(C)) \to \caH(C)$$
induced by the inclusion $[0] \to \Ar([2])$, $0 \mapsto \id_2$.
Then the fiber of the natural map $$F^\sim \to (C_{/X})^\sim$$
over a map $f \colon Y \to X$ is naturally equivalent to the subspace of the space
of paths in $\map_C(Y,Y^\vee)^{hC_2}$ from the zero map to the map
$$Y \overset{f}{\to} X \simeq X^\vee \overset{f^\vee}{\to} Y^\vee$$
on those connected components which exhibit the resulting commutative
square
$$\xymatrix{Y \ar[r] \ar[d] & X \ar[d] \\ 0 \ar[r] & Y^\vee}$$
as an exact square in $C$.
\end{proposition}

\begin{proof}
This follows from the above considerations together with the next two Lemmas.
\end{proof}

\begin{lemma}
\label{jttzzth}
Let $D$ be the category associated to the $1$-skeleton of $[1] \times [1]$ (so $D$ is Joyal-equivalent to two $\Lambda_1^2$'s glued together along
their endpoints) and $i \colon D \to [1] \times [1]$ the natural map. Then for an $\infty$-category $C$
the fiber of the functor $$i^* \colon \Fun([1] \times [1], C) \to \Fun(D,C)$$
over an object $\alpha \colon D \to C$ in the target is naturally equivalent to the space of paths in $\map(X,Y)$
($X$ being $\alpha((0,0))$ and $Y$ being $\alpha((1,1))$) from the composition of one composable pair of
maps in $D$ to the composition of the other composable pair.
\end{lemma}

\begin{proof}
This follows from the fact that there is a pushout square
$$\xymatrix{E \ar[r] \ar[d] & D \ar[d] \\ [1] \ar[r] & [1] \times [1],}$$
where $E$ is obtained by gluing two copies of $[1]$ together along their endpoints.
\end{proof}

\begin{lemma}
Let $C$ be an $\infty$-category which has a zero object.
Let $i \colon [2] \to [1] \times [1]$ be the map which sends $0$ to $(0,0)$, $1$ to $(0,1)$ and $2$ to $(1,1)$.
Let $$\Fun'([1] \times [1],C) \subset \Fun([1] \times [1],C)$$
be the full subcategory on those squares such that the entry in spot $(1,0)$ is a zero object.
Then the fiber of the functor $$\Fun'([1] \times [1],C) \to \Fun([2],C)$$
over an object $\alpha \colon [2] \to C$ is naturally equivalent to the space of paths
in the mapping space $\map(\alpha(0),\alpha(2))$ from the zero map to $\alpha(0 \to 2)$.

If $C$ has a duality and $\alpha \in \caH(\Fun([2],C))$, then
$\map(\alpha(0),\alpha(2))$ has a natural $C_2$-action, the map $\alpha(0 \to 2)$
naturally lies in $\map(\alpha(0),\alpha(2))^{hC_2}$ and the fiber
of the functor $$\caH(\Fun'([1] \times [1],C) \to \caH(\Fun([2],C))$$
is naturally equivalent to the space of paths in $\map(\alpha(0),\alpha(2))^{hC_2}$
from the zero map to $\alpha(0 \to 2)$.
\end{lemma}

\begin{proof}
The first part follows from Lemma \ref{jttzzth}, the second by taking hermitian objects.
\end{proof}

\begin{corollary}
Let $C \in (\Cat_\infty^{hC_2})^\st$ and $X \in \caH(C)$.
Let $F$ be the fiber over $X$ of the functor $$\caH(S_2(C)) \to \caH(C)$$
induced by the inclusion $[0] \to \Ar([2])$, $0 \mapsto \id_2$.
Then giving a point in $F$ is the same as giving a Lagrangian of $X$ in the sense of \cite[Example 7.]{lurie.L5}.
\end{corollary}

\begin{corollary}
\label{ht5z54t}
Let $C \in (\Cat_\infty^{hC_2})^\st$, $X \in \caH(C)$ and $\varphi \colon X \to X^\vee$ the corresponding
map in $$\map(X,X^\vee)^{hC_2} \simeq \Omega^\infty \map^\Sp(X,X^\vee)^{hC_2}.$$
Choose an inverse $-\varphi$ of $\varphi$ with respect to the infinite loop space structure
on $\map(X,X^\vee)^{hC_2}$. Then there is an object of $\caH(S_2(C))$
whose underlying exact triangle in $C$ has the form
$$\xymatrix{X \ar[r]^{\mathrm{diag}} \ar[d] & X \oplus X \ar[d]^{\varphi + (-\varphi)} \\
0 \ar[r] & X^\vee}$$
and where the hermitian structure on $X \oplus X$ is given by $\varphi \oplus (-\varphi)$.
\end{corollary}

\begin{proof}
This follows now from the proof of \cite[Proposition 11.]{lurie.L5}.
\end{proof}

\begin{proposition}
\label{wet5u6}
For $C \in (\Cat_\infty^{hC_2})^\st$ the $E_\infty$-structure on $\GW(C)$ defined in
section \ref{jtez4t4x} is grouplike.
\end{proposition}

\begin{proof}
There is a coequalizer diagram
$$\pi_0(S_1^{e,\sim}(C)_h) \rightrightarrows \pi_0(S_0^{e,\sim}(C)_h) \to \pi_0(|S_\bullet^{e,\sim}(C)_h|)$$
in $\Set$.
Let a point in $\pi_0(S_0^{e,\sim}(C)_h)$ be represented by a hermitian object $X \in \caH(C)$.
Let $\varphi \in \map(X,X^\vee)^{hC_2}$ be the corresponding map with a choice of an inverse
$-\varphi$. Let $W' \in \caH(S_2(C))$ be the object described in Corollary \ref{ht5z54t}
and $W$ be the image of $W'$ under the functor $\caH(S_2(C)) \to \caH(S_3(C))$ induced
by the map $[3] \to [2]$ given by $0 \mapsto 0, 1 \mapsto 1, 2 \mapsto 1, 3 \mapsto 2$.
The object $W$ determines an element of $\pi_0(S_1^{e,\sim}(C)_h)$ which is sent under the two maps above
to $0$ resp. $(X \oplus X, \varphi \oplus (-\varphi))$. So we see that in $\pi_0(|S_\bullet^{e,\sim}(C)_h|)$
an inverse of the image of $(X,\varphi)$ is given by $(X,-\varphi)$, in particular the $E_\infty$-space
$|S_\bullet^{e,\sim}(C)_h|$ is grouplike. It follows that also $\GW(C)$ is grouplike.
\end{proof}

\begin{definition}
For $n \in \naturals$ the abelian group $\GW_n(C):=\pi_n \GW(C)$ is called the $n$-th Grothendieck-Witt group
of the stable $\infty$-category $C$ with duality. In particular the group $\GW_0(C)$ is called
the Grothendieck-Witt group of $C$.
\end{definition}

\section{Hyperbolic categories}

We denote a right adjoint of the forgetful functor $$\Cat_\infty^{hC_2} \to \Cat_\infty$$
by $\Hyp$ and call $\Hyp(C)$ the hyperbolic category associated to $C \in \Cat_\infty$.
The underlying category of $\Hyp(C)$ is equivalent to $C \times C^\op$, and the duality is
informally given by $C \times C^\op \ni (X,Y) \mapsto (Y,X)$.

\begin{lemma}
Let $C \in \Cat_\infty$. Then there is a natural equivalence $$\caH^\lax(\Hyp(C)) \simeq \Tw(C)$$
of $\infty$-categories.
\end{lemma}

\begin{proof}
The $\infty$-category $\caH^\lax(\Hyp(C))$ is given as the complete Segal space
$$[n] \mapsto \map_{\Cat_\infty^{hC_2}}([n] * [n]^\op,\Hyp(C)),$$
which by adjunction is equivalent to $$[n] \mapsto  \map_{\Cat_\infty}([n] * [n]^\op,C).$$
But this is a possible definition of the twisted arrow category $\Tw(C)$.
\end{proof}

\begin{corollary}
There is a natural equivalence $$\caH(\Hyp(C))^\sim \simeq C^\sim$$
in $\Spc$ for $C \in \Cat_\infty$.
\end{corollary}

\begin{proof}
The core groupoid of the full subcategory of $\Tw(C)$ on the equivalences is naturally equivalent to $C^\sim$.
\end{proof}

\begin{lemma}
There is a natural equivalence $$\Fun(I,\Hyp(C)) \simeq \Hyp(\Fun(I,C))$$
in $\Cat_\infty^{hC_2}$ for $I \in \Cat_\infty^{hC_2}$ and $C \in \Cat_\infty$.
\end{lemma}

\begin{proof}
By adjunction a map $$\Fun(I,\Hyp(C)) \to \Hyp(\Fun(I,C))$$ is the same as a map
$$\Fun(I,\Hyp(C)) \to \Fun(I,C)$$ in $\Cat_\infty$, and such a map is induced
by the counit $\Hyp(C) \to C$. One checks that the resulting map in $\Cat_\infty^{hC_2}$
is an equivalence.
\end{proof}

\begin{corollary}
There is a natural equivalence $$\caH(\Fun(I,\Hyp(C)))^\sim \simeq \Fun(I,C)^\sim$$
in $\Cat_\infty$ for $I \in \Cat_\infty^{hC_2}$ and $C \in \Cat_\infty$.
\end{corollary}

It follows

\begin{lemma}
\label{htt53z5}
For $C \in \Cat_\infty^\st$ the simplicial object $S_\bullet^{e,\sim}(\Hyp(C))_h$ in $\Spc$ is naturally
equivalent to the simplicial object $S_\bullet^{e,\sim}(C)$.
\end{lemma}

\begin{proposition}
\label{j4t43wef}
For $C \in \Cat_\infty^\st$ there is a natural equivalence $\GW(\Hyp(C)) \simeq K(C)$
of grouplike $E_\infty$-spaces.
\end{proposition}

\begin{proof}
Note first that for any simplicial space $X$ the natural map $E(X) \to X$ induces an equivalence
$|E(X)| \to |X|$ (this follows from \cite[Lemma 1]{schlichting.mv}).
Thus by Lemma \ref{htt53z5} $\GW(\Hyp(C))$ is given as the fiber of the map
$$|S_\bullet^{e,\sim}(C)| \to |S_\bullet^{e,\sim}(C)| \times |S_\bullet^{w,\sim}(C^\op)| \simeq |S_\bullet^{e,\sim}(C)| \times |S_\bullet^{e,\sim}(C)|$$
which is naturally equivalent to the diagonal. The claim follows.
\end{proof}

\section{A connective real $K$-theory spectrum}

Recall from \cite[\S 6]{hls} the left adjoints $$\Cat \overset{f}{\to} \Cat^\Sigma \overset{l}{\to} \Cat^\preadd,$$ where $l$
is a localization. The functor $f$ sends a small $\infty$-category $C$ to the full subcategory of $\caP(C)$ which contains
the essential image of $C$ under the Yoneda embedding $C \to \caP(C)$ and is closed under finite coproducts,
see the proof of \cite[Proposition 5.3.6.2]{lurie.HTT}.


For $C \in \Cat_\infty^\Sigma$ we let $\caP_\Sigma(C):= \Fun^{\prod}(C^\op,\Spc)$, see \cite[\S 5.5.8]{lurie.HTT}.
For $C \in \Cat_\infty$ we let $$\caP^\oplus(C):= \Fun(C^\op,E_\infty(\Spc)) \simeq E_\infty(\caP(C))$$
(for the last equivalence see \cite[Remark 2.1.3.4]{lurie.higheralgebra}),
and for
$C \in \Cat_\infty^\Sigma$ we set $$\caP_\Sigma^\oplus(C):= \Fun^{\prod}(C^\op,E_\infty(\Spc)) \simeq E_\infty(\caP_\Sigma(C)).$$
If $C \in \Cat_\infty^\preadd$ then $\caP_\Sigma^\oplus(C) \simeq \caP_\Sigma(C)$, see \cite[Corollary 2.5 (iii)]{gepner-groth-nikolaus}.

For $C \in \Cat_\infty$ we have a canonical left adjoint $\caP(C) \to \caP^\oplus(C)$.

For $C \in \Cat_\infty^\Sigma$ we have a natural square of left adjoints
$$\xymatrix{\caP(C) \ar[r] \ar[d] & \caP_\Sigma(C) \ar[d] \\
\caP^\oplus(C) \ar[r] & \caP_\Sigma^\oplus(C)}$$
which commutes up to a natural equivalence since the corresponding right adjoints do.
In particular we exhibit a natural functor $C \to \caP_\Sigma^\oplus(C)$.

For $C \in \Cat_\infty$ let $F^\oplus(C)$ be the smallest full subcategory of $\caP^\oplus(C)$ that contains the essential
image of $C \to \caP^\oplus(C)$ and is closed under finite coproducts. Note $F^\oplus(C)$ is preadditive.
Similarly for $C \in \Cat_\infty^\Sigma$ let
$F_\Sigma^\oplus(C)$ be the essential image of $C \to \caP_\Sigma^\oplus(C)$. $F_\Sigma^\oplus(C)$ is also preadditive.

\begin{proposition}
i) There is a natural equivalence of functors $l \circ f \simeq F^\oplus$, so for $C \in \Cat_\infty$
the category $F^\oplus(C)$ is the free preadditive category on the $\infty$-category $C$.

ii) There is a natural equivalence of functors $l \simeq F_\Sigma^\oplus$, so for $C \in \Cat_\infty^\Sigma$
the category $F_\Sigma^\oplus(C)$ is the free preadditive category on the $\infty$-category $C$ with finite coproducts.
\end{proposition}

\begin{proof}
The first point is \cite[Proposition 2.8]{glasman.good}, the second point follows similarly.
\end{proof}

For a small $\infty$-category $C$ which has pullbacks we denote by $\Span(C)$
the $\infty$-category of spans in $C$. That is the $\infty$-category denoted $A^{\mathit{eff}}(C)$ in
\cite[Definition 3.6]{barwick.mack1}.

The assignment $C \mapsto \Span(C)$ can be viewed as a functor
$$\Span \colon \Cat_\infty^\lex \to \Cat_\infty,$$
where $\Cat_\infty^\lex$ is the subcategory of $\Cat_\infty$ of $\infty$-categories with all finite limits and left exact functors
between them (see loc. cit.).

We denote by $\Fin$ the $\infty$-category of finite sets.

\begin{proposition}
\label{h5r34ef}
The $\infty$-category $\Span(\Fin)$ is preadditive, and the natural functor
$F^\oplus(*) \to \Span(\Fin)$ sending the point to the one element set is an equivalence.
\end{proposition}

\begin{proof}
This follows by comparing mapping spaces.
\end{proof}

We recall the

\begin{theorem}
\label{hte45z5}
Let $G$ be a finite group. Then the symmetric monoidal $\infty$-category of genuine $G$-spectra is equivalent to the full subcategory
$$\Fun^{\coprod}(\Span(\Fin[G]),\Sp) \subset \Fun(\Span(\Fin[G]),\Sp)$$ of finite coproduct preserving functors
equiped with the symmetric monoidal structure which is induced by the Day convolution product on the functor category.
\end{theorem}

\begin{proof}
The equivalence is \cite{guillou-may}, \cite[Example B.6]{barwick.mack1}, the symmetric monoidal structure is (partially) discussed
in \cite[\S 3]{barwick.mack2}.
\end{proof}

\begin{lemma}
\label{j6tr457yu}
Let $C$ be a cocomplete $\infty$-category. Then the functor $F \colon C \to \Fun(\Spc,C)$ which sends $X \in C$ to
the functor $\Spc \ni K \mapsto K \otimes C$ is left adjoint to the functor $\Fun(\Spc,C) \ni \varphi \mapsto \varphi(*)$.
\end{lemma}

\begin{proof}
This follows from the fact that for any $X \in C$ the functor $F(X)$ is the left Kan extension of
the functor $* \to C$ which sends the unique object of $*$ to $X$ along the inclusion $* \to \Spc$ which sends the object of $*$
to $*$.
\end{proof}

Let $C \in \Cat_\infty^\lex$ and $\varphi \colon \Spc \to \Cat_\infty^\op$ the functor $K \mapsto \Span(\Fun(K,C))$
(note that $K \mapsto \Fun(K,C)$ is a functor $\Spc \to (\Cat_\infty^\lex)^\op$).
Then $\Span(C) \simeq \varphi(*)$, and this equivalence is by Lemma \ref{j6tr457yu} adjoint to a map
$F \to \varphi$ in $\Fun(\Spc,\Cat_\infty^\op)$, where $F$ is defined by $\Spc \ni K \mapsto \Fun(K,\Span(C))$.

In particular this exhibits for any $K \in \Spc$ a natural functor $$f_C^K \colon \Span(\Fun(K,C)) \to \Fun(K,\Span(C)).$$

For a finite group $G$ we set $g_G:= f_\Fin^{BG}$, so we have $$g_G \colon \Span(\Fin[G]) \to \Span(\Fin)[G].$$

\begin{lemma}
The $\infty$-category $(\Cat_\infty^{hC_2})^\st$ is preadditive.
\end{lemma}

So by Proposition \ref{h5r34ef} for $C \in (\Cat_\infty^{hC_2})^\st$ we exhibit a natural
functor $$\Span(\Fin) \to (\Cat_\infty^{hC_2})^\st$$
sending the generator to $C$,
which induces the second map in the composition
$$\Span(\Fin[C_2]) \overset{g_{C_2}}{\longrightarrow} \Span(\Fin)[C_2] \to (\Cat_\infty^{hC_2})^\st[C_2]
\to  (\Cat_\infty^{hC_2})^\st \overset{\GW}{\longrightarrow} \Sp,$$
whereas the third map is the map $t_{\Cat_\infty^\st}$ described in \cite[\S 9]{hls} (and $\GW$ takes
values in grouplike $E_\infty$-spaces aka connective spectra).

By Theorem \ref{hte45z5} we obtain a genuine $C_2$-spectrum which we denote by $\KR(C)$ and call
the connective real $K$-theory spectrum of $C$.

\begin{proposition}
For $C \in (\Cat_\infty^{hC_2})^\st$ the underlying spectrum of $\KR(C)$ is naturally equivalent to
the $K$-theory spectrum $K(C)$ (which thus inherits a natural $C_2$-action), and the $C_2$-fixed points
of $\KR(C)$ are canonically equivalent to $\GW(C)$.
\end{proposition}

\begin{proof}
In the defining composition of $\KR(C)$ the $C_2$-set $*$ is sent to the $C_2$-fixed points.
But the image of $*$ in $(\Cat_\infty^{hC_2})^\st[C_2]$ is $C$ with the trivial $C_2$-action
yielding the second claim.

The underlying spectrum of $\KR(C)$ is the image of the $C_2$-orbit $C_2$. Its image
in $(\Cat_\infty^{hC_2})^\st[C_2]$ is $C \times C$ with the $C_2$-action which switches the two factors. The resulting object
$t_{\Cat_\infty^\st}(C \times C)$ of $(\Cat_\infty^{hC_2})^\st$ is then seen to be naturally equivalent
to $\Hyp(C)$, so the first claim follows from Proposition \ref{j4t43wef}.
\end{proof}

\bibliographystyle{plain}
\bibliography{gw}

\vspace{0.1in}

\begin{center}
Fakult{\"a}t f{\"u}r Mathematik, Universit{\"a}t Osnabr\"uck, Germany.\\
e-mail:\\
markus.spitzweck@uni-osnabrueck.de
\end{center}

\end{document}